# A sustainable inventory model considering a discontinuous transportation cost function and different sources of pollution

**Alfonso Angel Medina-Santana**

Department of Industrial and Systems Engineering, School of Engineering and Sciences, Tecnológico de Monterrey, E. Garza Sada 2501 Sur, Monterrey, Nuevo León, CP 64 849, México

A01670119@itesm.mx

**Leopoldo Eduardo Cárdenas-Barrón**

Department of Industrial and Systems Engineering, School of Engineering and Sciences, Tecnológico de Monterrey, E. Garza Sada 2501 Sur, Monterrey, Nuevo León, CP 64 849, México

lecarden@tec.mx

**Abstract**

The economic ordering quantity (EOQ) is the first inventory model which has still being studying extensively till nowadays. Several considerations have been incorporated throughout the time such as backorders and imperfect quality. More recently, sustainability issues have also been introduced to the literature of inventory models remarking the increasing importance of the three bottom line: environment, social and economic criteria. In that sense, the aim of this research work is to develop a new sustainable inventory model considering a discontinuous transportation cost function and a new carbon emissions function to include different sources of pollution in the decision making process that are not taken into account using the traditional approaches in the inventory literature. The proposed inventory model permits to obtain a continuous and a discrete solution for the lot size and constitutes an easy-to-use process that can be implemented in a simple manner.

**Keywords:** Sustainable EOQ; emissions function; discontinuous transportation cost function; continuous and integer lot size.

# 1. Introduction

The concerns and efforts on reducing greenhouse gases (GHG) emissions and their corresponding effect on global warming have incremented through the years. The Kyoto protocol was established in the United Nations framework in 1997 and it was ratified by 156 countries. According to Almer and Winkler (2017) the main objective of Kyoto protocol is to establish the reduction of GHG emissions in 5.2% with respect to 1990. Later, in 2015, the Paris agreement was introduced and entered into vigor in 2016. This protocol was approved by 195 countries and its aim is to keep the increment of global temperature below 2 degrees Celsius and to pursue efforts to limit the temperature rise below to 1.5 degrees Celsius (United Nations-climate change, 2018). Additionally, it seeks to strengthen the ability of countries to deal with the impacts of climate change. In that sense, each country makes efforts to reduce the environmental impact in some specific areas. For example, Wang, et al. (2019) states to consider environmental issues in logistic activities. Since Harris, (1913) formulated the economic order quantity (EOQ) inventory model in 1913, different characteristics, including sustainability criteria, have been incorporated into a common framework which is the EOQ inventory model. Arslan and Turkay (2013) revised the standard EOQ inventory model to integrate sustainability considerations that take account of environmental issues. In that work the authors include objectives and constraints; and they analyze five different approaches: Direct accounting, carbon tax, direct cap, cap and trade and carbon offsets. As it is depicted in the latter work, the most practical approach is the direct accounting approach because the sustainability issues are translated into costs and modeled as a part of the total cost function.

Other approaches are useful to study governmental policies. For example, Hua et al. (2011a) studied a sustainable inventory model under the cap-and-trade mechanism, comparing the classical EOQ inventory model with a sustainable one. They examine the impact of carbon cap and carbon price on order size as a method of inducing the retailer to reduce carbon emissions that result in the increase of total costs. Bonney and Jaber (2011) discussed in depth how to design responsible inventory systems including not only the importance of some activities in the total cost of the system but also their effect in the environment and interrelations with other activities. Additionally, a simplistic model that extends the EOQ to include some environmental costs is proposed. Chen et al. (2013) analyzed if it is possible to reduce emissions by changing the order quantities. Moreover, they provide a condition under which the reduction in carbon emissions is greater than the increment in cost and studied which factors affect the difference between these two quantities.

The most common sustainability criterion is the amount of carbon emissions and it is included in several works as a cost function or a constraint. Absi et al. (2013) proposed carbon emission constraints in a multi-sourcing lot-sizing problem. The aim of these constraints is to put a limit in the carbon emission per unit of product supplied. Furthermore, this inventory model considers a multi-sourcing lot-sizing problem in which a company must determine, over a planning horizon of $T$ periods, when, where and how much to produce of an item in order to satisfy a deterministic time-dependent demand. Moreover, different production locations and transportation modes are available to satisfy the demand. Benjaafar et al. (2013) proposed an inventory model that incorporates carbon emission constraints on single and multi-stage lot-sizing models with a cost minimization objective. For this purpose, four regulatory policy settings are considered, based respectively on a strict carbon cap, a tax on the amount of emissions, the cap-and-trade system and the possibility to invest in carbon offsets to mitigate carbon caps.

Other studies have focused on the optimization of the inventory management and the supply chain. Saadany et al. (2011) and Jaber et al. (2013) proposed formulations in which emissions from manufacturing processes are included into a two-echelon supply chain model. The first study assumes that demand is a function of the price and product's environmental quality whereas the second one discusses the efficiency of different emissions trading schemes by numerical examples. Wahab et al. (2011) found the optimal production-shipment policies by proposing three scenarios for a two-level supply chain. They also include the environmental impact by considering the fixed and variable carbon emission costs. In this work, the objective is to minimize the expected total cost per unit of time.

A large number of authors explore approaches different from the EOQ. For instance, Rosič and Jammernegg (2013) extended the dual sourcing model based on the newsvendor framework by considering the environmental impact of transport. In the first step, the effect of a linear carbon emission tax on transport is analyzed, and in the second step, an emission trading scheme for the transport sector is proposed. Choi and Chiu (2012) explored the mean-downside-risk (MDR) and mean-variance (MV) newsvendor models under both the exogenous and endogenous retail price decision cases. First, analytical models with the MDR and MV objectives are constructed. Then, it is shown that the analytical solution schemes for both the MDR and MV problems are the same. In contrast, a different approach is proposed by Glock (2012) who focuses on coordinated inventory replenishment decisions between buyer and vendor and their impact on the performance of the supply chain. The joint economic lot size (JELS) models determine simultaneously the order, production and shipment quantities from the perspective of the supply chain with the objective of minimizing total system costs whereas carbon emissions constraints are considered.

Several formulations and optimization techniques have been implemented to solve certain sustainable inventory models. Battini et al. (2014) included a discontinuous transportation cost function and the corresponding optimization algorithm to solve an inventory problem with discontinuity points. Furthermore, the environmental impact of transportation and inventory management are incorporated into the model and investigated from an economic point of view. Additionally, internal and external transportation costs, vendor and supplier location and the different freight vehicle utilization ratio are considered to include the transportation mode selection. On the other hand, a multi objective model is presented by Bouchery et al. (2012). In their work, it is underlined that reducing all aspects of sustainable development to a single objective is not desirable and reformulate the classical economic order quantity model EOQ as a multi objective problem. Finally, to the best of our knowledge, there exists just one multi-product, multi-echelon sustainable inventory model that is shown by Zhang and Xu (2013).

More recently, plenty of sustainable inventory models have been developed to fill some gaps in the literature. For instance, Zadjafar and Gholamian (2018) included environmental ergonomics based on the argument which states that the three pillars of sustainability (economical, environmental and social) are not independent among them but interdependent. The proposed methodology is useful to include the effects of the environment on the social aspect, more specifically human health because of the emission of different gases. Other example is the work proposed by Taleizadeh et al. (2018) which developed four economic production quantity (EPQ) models which consider different shortage situations by using the direct accounting approach. From the numerical experiments of this study, it is concluded that the partial backordering situation has enough generality and usability for real business environments. Wang et al. (2019) studied the inventory and pricing policies from an environmental perspective. The authors considered a single period inventory model with recycling action and the government subsidy for one and two supply option systems. Then, the effects of environmentally sensitive for customers and the green degree of products on the pricing policies are studied. Finally, Tiwari et al. (2018) included deteriorating and imperfect quality items to a sustainable inventory model which seeks to minimize a total cost function which integrates the buyer's and the vendor's costs.

In this research work, a new sustainable inventory model which includes a discontinuous transportation cost function and different sources of pollution is formulated and solved. It is important to remark that this inventory model also incorporates a new component in the carbon emissions function. To optimize the total cost function, a Taylor polynomial is used and therefore, the optimization procedure is based on an approximation. This approximation is analyzed, and a sensitivity analysis is carried out to figure out how the accuracy of the approximation changes as the value of parameters increases or decreases. Additionally, an algorithm and a numerical experiment are carried out when the order quantity is considered as a continuous and discrete variable. From this mathematical procedure, it is worth highlighting that a useful approximation is obtained which not only allows to minimize the cost of the inventory system when the lot size is a continuous variable but also to optimize the inventory policy when the lot size is an integer value. Moreover, by using this approximation the proposed optimization procedure is easy to use and reduce considerably the computational complexity.

The structure of this research work is as follows: Section 2 presents the notation and assumptions. Section 3 formulates the mathematical inventory model which comprehends environmental and non-environmental costs, and identifies the special cases to which the proposed model reduces. Section 4 develops the continuous and discrete optimization algorithms. Section 5 solves some numerical examples. Finally, Section 6 provides some conclusions and future research work directions.

## 2. Notation and assumptions

This section contains the notation and assumptions that are utilized in the development of the sustainable inventory model.

### 2.1. Notation

The following notation is as follows:

$A$ : Ordering cost ($/order)

$c$ : Unit purchasing cost (\$/unit)

$h$ : Holding cost (\$/unit/units of time)

$a$ : Fixed cost per trip (\$/order)

$b$ : Variable cost per unit transported per distance travelled (\$/unit/units of distance)

$d$ : Distance traveled from supplier to buyer (units of distance/order)

$\alpha$ : Proportion of waste returned ($0 \leq \alpha \leq 1$)

$D$ : Demand rate (units/units of time)

$Q$ : Lot size (units/order)

$\beta$ : Social cost from vehicle emission (\$/units of time)

$v$ : Average velocity (units of distance/units of time)

$\gamma$ : Cost of disposing waste to the environment (\$/unit)

$\gamma_0$ : Fixed cost of disposing waste to the environment (\$/order)

$\theta$ : Proportion of waste produced per lot $Q$ ($0 \leq \theta \leq 1$)

$\varepsilon$ : Units of carbon emissions per replenishment (units of $CO_2$/order)

$g$ : Unit variable emissions in the warehouse (units of $CO_2$/ unit)

$C_e$ : Cost incurred due to the emissions when ordering and holding (\$/units of $CO_2$)

$S$ : Saturation level of the container ($0 \leq S \leq 1$)

$i$ : Type of container

$n_i$ : Number of containers type $i$

$y_i$ : Capacity of the container type $i$

$C_p$ : Cost per unit of capacity of the container (\$/units of capacity)

$l$ : Shape parameters for the emissions curve of the industry in case (units of $CO_2$/unit/units of time).

$r$ : Shape parameters for the emissions curve of the industry in case (units of time/order)

A useful explanation of the shape parameters $r$ and $l$ is given in appendix A.

## 2.2. Assumptions

The sustainable inventory model presented in this study has the following assumptions:

a) The demand is known, constant and continuous.
b) The lead time is known and constant.
c) The replenishment is instantaneous.

Following these assumptions, the formulation of the mathematical model is carried out in Section 3.

## 3. The mathematical formulation of the sustainable inventory model

The sustainable inventory model includes environmental and traditional costs into a total cost function. Both costs are described in detailed manner below.

### 3.1. Environmental costs

According to Aronsson and Brodin (2006), the measurement of emissions is one of the most important indicators for the environmental impact, for this reason the carbon emissions are considered in this paper. Commonly, the formula given in equation (1) is used to calculate the carbon emissions during a period (e.g. a year).

$$CO_2(Q) = \varepsilon \frac{D}{Q} + g \frac{Q}{2} \tag{1}$$

This formula is widely used in the literature; see for example the research works of Bouchery et al. (2012), Hovelaque and Bironneau (2015), and Hua et al. (2011b). As it is expressed in equation (1), the inventory system incurs in an amount of emissions each time there is a replenishment and also incurs in an amount of emissions related to the holding stock. The emissions ($\varepsilon$) per replenishment are associated with the carbon emissions of activities done during each cycle such as the use of some machines. An amount $g$ of emissions per unit of average inventory level is incurred too. This is due to the fact that a greater amount of stock requires more illumination, air conditioning, continue use of some machines to manipulate the inventory. It is evident that all these activities incur in an amount of carbon emissions.

With the aim of accounting the carbon emissions more closely to real life situations, this work proposes a new formula which is shown in equation (2).

$$CO_2(Q) = \varepsilon \frac{D}{Q} + g \frac{Q}{2} + l \frac{Q}{2} e^{rD/Q} \tag{2}$$

The difference is in the new expression $l\frac{Q}{2}e^{rD/Q}$ which suggests that one of the sources of emissions during a period is affected by the average level of inventory proportionally and the number of replenishments exponentially. This new expression corresponds to an amount of emissions that cannot be determined using the formula in equation (1). It is considered then the expression $l\frac{Q}{2}e^{rD/Q}$ as a surplus of emissions. The Figure 1 helps to illustrate the benefits of using this expression.

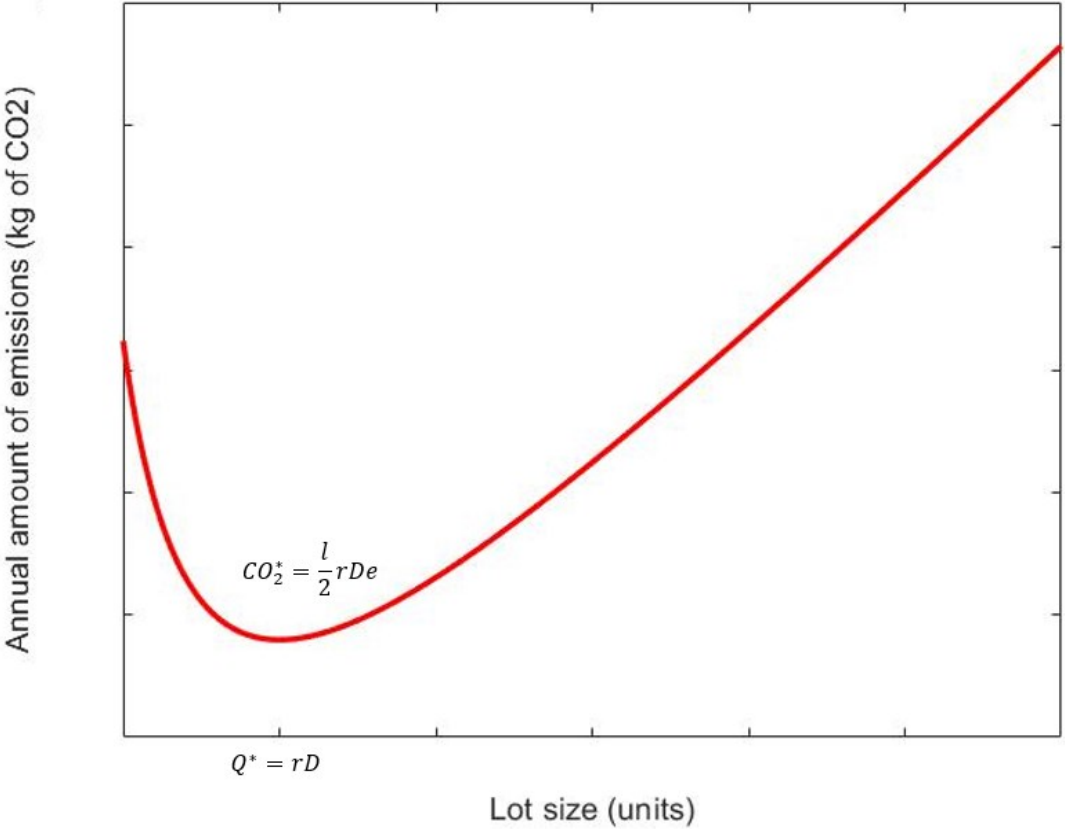

**Figure 1.** Emissions calculated by using equation (2)

It is appreciated in the graph that, for a lot size greater than $Q^* = rD$, the amount of emissions during a period tends to be proportional to the parameter $l$. However, if the lot size is lower than this critical value, the amount of emissions starts increasing rapidly in the direction in which the lot size is diminishing. This can be an attribute of the carbon emissions behavior of some machines whose emissions are related to their energy consumption (that is proportional to the average inventory level during a period). However, if it is assumed that these machines are turned off and turned on every cycle, the corresponding carbon emissions can grow rapidly since a small batch size implies many orders over a period which means, for instance, to turn off and turn on the air conditioner several times diminishing its performance.

Dividing equation (2) by $D/Q$, the emissions per cycle are given as:

$$CO_{2p.c}(Q) = \varepsilon + \frac{gQ^2}{2D} + \frac{lQ^2}{2D} e^{rD/Q} \tag{3}$$

Then, the related cost with the amount of these emissions is therefore expressed:

$$C_1(Q) = CO_{2p.c}(Q)C_e \tag{4}$$

$$C_1(Q) = \left(\varepsilon + \frac{gQ^2}{2D} + \frac{lQ^2}{2D} e^{rD/Q}\right)C_e \tag{5}$$

Additionally, there is also a part of emissions related to the transportation process. Thus, since the emissions are directly related to the time that vehicle is traveling, thus the emissions cost from transportation is formulated as in Bonney and Jaber (2011).

$$C_2(Q) = 2\beta \frac{d}{v} \tag{6}$$

Likewise, the waste disposal process has also an environmental impact. It is assumed that there is an amount of waste produced by the company during each cycle; this amount is $Q\theta$, and the customer uses the company's transportation system to return an amount of waste $Q\alpha$ each time that the units are supplied in order to meet the demand. For convenience, both quantities are considered as an independent proportion of the lot size. Besides, a fixed cost $\gamma_0$ and a variable cost $\gamma$ are incurred because of the disposal activity. This is expressed as follows:

$$C_3(Q) = \gamma_0 + \gamma Q\theta + \gamma Q\alpha \tag{7}$$

$$C_3(Q) = \gamma_0 + \gamma Q(\theta + \alpha) \tag{8}$$

The use of containers is also considered for the environmental impact. In that sense, it is supposed that the index $i$ corresponds to a certain type of container and each type has a capacity of $y_i$ and a quantity $n_i$ of them is used. Then, the saturation level of each delivery is:

$$S = \frac{Q}{\sum n_i y_i} \tag{9}$$

Furthermore, the cost of the containers is related to the saturation level and the lot size.

$$C_4(Q) = C_p \frac{Q}{S} \tag{10}$$

$$C_4(Q) = C_p \sum_i n_i y_i \tag{11}$$

$$C_4(Q) = C_p D P_j \tag{12}$$

Since the set $N = \{n_i\}$ varies indicating that the customer can select among different combinations of the number of containers from each type, each possible combination represents a discontinuity point $DP_j$ in which the function is not differentiable since it is not continuous. Then, the optimization procedure requires a step by step approach as in Battini et al. (2014).

### 3.2. Non-environmental costs

Traditional costs of ordering $A$, holding $h$ and purchasing $c$ are considered as in the EOQ basic inventory model proposed by Harris (1913).

$$C_5(Q) = A + cQ + h\frac{Q^2}{2D} \tag{13}$$

Finally, the transportation cost related to the delivery cycle that includes the collection of the returned waste is expressed as follows:

$$C_6(Q) = 2a + bdQ + bd\alpha D \frac{Q}{D} \tag{14}$$

$$C_6(Q) = 2a + bdQ + bd\alpha Q \tag{15}$$

$$C_6(Q) = 2a + bdQ(1+\alpha) \tag{16}$$

### 3.3. Total cost function

Then, the total cost function per unit of time is:

$$M(Q) = \frac{1}{T} \sum_{i=1}^{6} C_i \tag{17}$$

Where:

$$T = \frac{Q}{D} \tag{18}$$

Then, the total cost function is explicitly expressed as:

$$M(Q) = \frac{D}{Q}[(\varepsilon + \frac{gQ^2}{2D} + \frac{lQ^2}{2D} e^{rD/Q})C_e + 2\beta\frac{d}{v} + \gamma_0 + \gamma Q(\theta + \alpha) + C_P DP_k + A + cQ + h\frac{Q^2}{2D} + 2a + bdQ(1+\alpha)] \tag{19}$$

$$M(Q) = [(\varepsilon \frac{D}{Q} + g\frac{Q}{2} + l\frac{Q}{2}e^{rD/Q})C_e + 2\beta\frac{d}{v}\frac{D}{Q} + \gamma_0\frac{D}{Q} + \gamma Q(\theta + \alpha)\frac{D}{Q} + C_p DP_k\frac{D}{Q} + A\frac{D}{Q} + cQ\frac{D}{Q} + h\frac{Q}{2} + 2a\frac{D}{Q} + bdD(1+\alpha)] \quad (20)$$

Since the complex expression $l\frac{Q}{2}e^{rD/Q}$ does not permit to obtain a closed form to calculate the lot size, it is proposed the use of a Taylor polynomial to reduce it. Then, as it is previously known, the equation (21) is a Taylor series and using only the first three terms of this expression it is obtained a Taylor polynomial given in equation (22).

$$e^x = \sum_{i=0}^{\infty} \frac{x^i}{i!} \quad (21)$$

$$e^x \cong \frac{x^0}{0!} + \frac{x^1}{1!} + \frac{x^2}{2!} \quad (22)$$

Using this Taylor polynomial, it is reformulated the equation (2):

$$l\frac{Q}{2}e^{rD/Q} \cong l\frac{Q}{2}(1 + r\frac{D}{Q} + \frac{r^2 D^2}{2Q^2}) \quad (23)$$

$$l\frac{Q}{2}e^{rD/Q} \cong \frac{l}{2}(Q + rD + \frac{r^2 D^2}{2Q}) \quad (24)$$

Thus, by replacing this expression in the total cost function (20):

$$M(Q) = [(\varepsilon\frac{D}{Q} + g\frac{Q}{2} + \frac{l}{2}(Q + rD + \frac{r^2 D^2}{2Q}))C_e + 2\beta\frac{d}{v}\frac{D}{Q} + \gamma_0\frac{D}{Q} + \gamma Q(\theta + \alpha)\frac{D}{Q} + C_p DP_k\frac{D}{Q} + A\frac{D}{Q} + cQ\frac{D}{Q} + h\frac{Q}{2} + 2a\frac{D}{Q} + bdD(1+\alpha)] \quad (25)$$

The first derivate with respect to $Q$ is:

$$M'(Q) = -\frac{\varepsilon DC_e}{Q^2} + \frac{gC_e}{2} + \frac{lC_e}{2} - \frac{lr^2 D^2 C_e}{4Q^2} - \frac{2\beta dD}{vQ^2} - \gamma_0\frac{D}{Q^2} - \frac{C_p DP_k D}{Q^2} - \frac{AD}{Q^2} + \frac{h}{2} - 2\frac{aD}{Q^2} \quad (26)$$

The second derivative with respect to $Q$ is

$$M''(Q) = \frac{2\varepsilon DC_e}{Q^3} + \frac{lr^2 D^2 C_e}{2Q^3} + \frac{4\beta dD}{vQ^3} + \frac{2\gamma_0 D}{Q^3} + \frac{2C_p DP_k D}{Q^3} + \frac{2AD}{Q^3} + \frac{4aD}{Q^3} > 0, \text{ for } Q > 0 \quad (27)$$

The order quantity is obtained easily by putting the first derivative equal to zero and solving for $Q$. Thus, the lot size is determined with

$$Q^* = \sqrt{\frac{2D(\varepsilon C_e + \frac{lr^2 DC_e}{4} + 2\beta\frac{d}{v} + \gamma_0 + C_p DP_k + A + 2a)}{gC_e + h + lC_e}} \quad (28)$$

Notice that the lot size $Q$ given by equation (28) minimizes the cost function in equation (25) if we consider it to be continuous in the open interval $(0,\infty)$ because in that case the convexity of the function would be proven by equation (27). Moreover, it must be remembered that the expression in equation (25) is an approximation of the total cost function which is expressed in equation (20). However, since both cost functions are discontinuous, the main procedure to obtain a continuous optimal solution involves the evaluation of the cost functions in a set of ranges. This algorithm is explained in Section 4.

### 3.4. Special cases

Considering only $C_5(Q)$ in equation (13) the proposed inventory model is reduced to Harris (1913)'s inventory model. Removing $C_1(Q)$ and $C_4(Q)$ the proposed inventory model reduces to Bonney and Jaber (2011)'s inventory model. Considering only $C_1(Q)$ and $C_5(Q)$ this work reduces to the direct accounting approach proposed in Arslan and Turkay (2013). For this a modification is required in equation (29)

$$C_1(Q) = (\varepsilon + gQ^2/2D + l(Q^2/2D)e^{(rD/Q)})C_e \qquad (29)$$

It is required to put $l = 0$. So, this expression is then equivalent to the following:

$$C_1(Q) = (\varepsilon + gQ^2/2D)C_e \qquad (30)$$

This makes that the proposed inventory model reduces to Arslan and Turkay (2013)'s model. In the following sections, the proposed formulation is used to face a discontinuous transportation cost function. Section 4 and Section 5 present algorithms for the case in which the lot size $Q$ is a continuous and an integer variable, respectively.

### 4. Case 1: when the lot size ($Q$) is a continuous variable

This section presents a main procedure and numerical results when the lot size is a continuous variable.

#### 4.1. Main procedure

This part includes three steps to find a continuous solution:

**Step 1.** Define several $j$ th ranges and for each range $j = \{DP_j, DP_{j+1}\}$, compute $Q'_j{}^*$ with equation (28).

**Step 2.** Then, in each range, the local optimal $Q_j{}^*$ is calculated as follows:

$Q_j{}^* = Q'_j{}^*$ when $DP_j < Q'_j{}^* \leq DP_{j+1}$

If $Q'_j{}^*$ is not into the range $j = \{DP_j, DP_{j+1}\}$, then

$Q_j^* = DP_j$ when $C(DP_j) \leq C(DP_{j+1})$

$Q_j^* = DP_{j+1}$ when $C(DP_j) > C(DP_{j+1})$

**Step 3.** At the end, the total minimum cost is defined as the minimum among the local total minimum cost values that satisfy the equation $C(Q^*) = \min\{C(Q_j^*)\} \forall j$ and so $Q^*$ is the optimal value for the lot size. In Subsection 4.2 a numerical example is carried out in order to clarify the use of the optimization algorithm.

### 4.2. Numerical experiments

With the aim to illustrate the use of the inventory model a numerical example is provided and solved. The data for the parameters of the example are given in Table 1.

**Table 1.** Parameters for the numerical experiment

| Parameter | Value |
|---|---|
| $A$ | 1000 ($/order) |
| $c$ | 25 ($/unit) |
| $h$ | 8 ($/unit/year) |
| $a$ | 80 ($/order) |
| $b$ | 4 ($/unit) |
| $d$ | 3000 (kilometers) |
| $\alpha$ | 0.1 |
| $D$ | 5000 (units/year) |
| $\beta$ | 30 ($/hour) |
| $v$ | 50 (kilometers/hour) |
| $\gamma$ | 5 ($/unit) |
| $\gamma_0$ | 20 ($/order) |
| $\theta$ | 0.1 |
| $\varepsilon$ | 200 (kg of $CO_2$/order) |
| $g$ | 3 (kg of $CO_2$/unit/year) |
| $C_e$ | 10 ($/kg of $CO_2$) |
| $C_p$ | 2 ($/unit of capacity) |
| $r$ | 0.004 (years/order) |
| $l$ | 30 (kg of $CO_2$/unit/year) |

Regarding the discontinuity points, it is considered that the supplier has two options of containers for transporting the product. The capacities of each container $(y_1, y_2)$ are (300,600). Moreover, since there are only two containers of each type available, six ranges are studied and enumerated in the first column in Table 2. Using the values of parameters in Table 1 and the algorithm (main procedure) previously described in Subsection 4.1, the obtained results are shown in Table 2 and Table 3. Table 2 is obtained by using the total cost function and in Table 3, the results are obtained using the proposed approximation total cost function. It must be observed that the results are pretty similar.

**Table 2**. Results using the total cost (TC) function

| $j+1$ | $n_1$ | $n_2$ | $DP_j$ | $DP_{j+1}$ | TC($Q'_j$*opt) | TC($Q'_j$*app) | TC($DP_j$) | TC($DP_{j+1}$) | $Q'_j$*opt | $Q'_j$*app | T/F | $Q_j$* |
|---|---|---|---|---|---|---|---|---|---|---|---|---|
| 1 | 1 | 0 | 0 | 300 | | | - | **66306802.260** | 467.468 | 467.462 | F | 300.0 |
| 2 | 0 | 1 | 300 | 600 | <u>**66297295.347**</u> | 66297295.347 | | | 486.084 | 486.078 | T | 486.0 |
| 3 | 1 | 1 | 600 | 900 | | | **66305950.560** | 66332800.248 | 504.012 | 504.008 | F | 600.0 |
| 4 | 2 | 1 | 900 | 1200 | | | **66336133.582** | 66374075.139 | 521.325 | 521.321 | F | 900.0 |
| 5 | 1 | 2 | 1200 | 1500 | | | **66376575.139** | 66419120.089 | 538.081 | 538.077 | F | 1200.0 |
| 6 | 2 | 2 | 1500 | 1800 | | | **66421120.089** | 66466050.062 | 554.331 | 554.327 | F | 1500.0 |

The main procedure is applied to obtain Table 2 and Table 3. For each range, a local minimum is obtained using equation (28) and it is allocated in column $Q'_j$*app for the $j$ th range. It must be reminded that this formula uses an approximation and this lot size is not the corresponding local minimum for the total cost function in (20). However, the true local optimal is obtained by using Lingo and is reported is column $Q'_j$*opt for the $j$ th range. Then, carrying out the optimization algorithm implies the evaluation of a cost function, for each range, in the local minimum found by using equation (28) (or Lingo) or the in the discontinuity points. It must be also realized that in case the local minimum is obtained using equation (28) or Lingo for a certain range, the condition true (T) or false (F) is reported in column T/F. There are two options to carry out the main procedure: to use the total cost function (TC) in equation (20) or the approximated total cost function (AC) in equation (25). Then, both options are implemented in Table 2 and Table 3, respectively.

**Table 3**. Results using the approximated cost function (AC)

| $j+1$ | $n_1$ | $n_2$ | $DP_j$ | $DP_{j+1}$ | AC($Q'_j$*app) | AC($DP_j$) | AC($DP_{j+1}$) | $Q'_j$*app | T/F | $Q_j$* |
|---|---|---|---|---|---|---|---|---|---|---|
| 1 | 1 | 0 | 0 | 300 | | - | **66306800.000** | 467.462 | F | 300.0 |
| 2 | 0 | 1 | 300 | 600 | <u>**66297294.492**</u> | | | 486.078 | T | 486.1 |
| 3 | 1 | 1 | 600 | 900 | | **66305950.000** | 66332800.000 | 504.008 | F | 600.0 |
| 4 | 2 | 1 | 900 | 1200 | | **66336133.333** | 66374075.139 | 521.321 | F | 900.0 |
| 5 | 1 | 2 | 1200 | 1500 | | **66376575.000** | 66419120.000 | 538.077 | F | 1200.0 |
| 6 | 2 | 2 | 1500 | 1800 | | **66421120.000** | 66466050.000 | 554.327 | F | 1500.0 |

In Figure 2, it is observed that there are several points of discontinuity (two in the graph). Moreover, it is proved that the optimal solution is $Q^*_{opt}$ = 486.084 from Table 2. By using the approximation, as it is proposed in the model formulation, the optimum is $Q^*_{app}$ = 486.078 which is reported in Table 3.

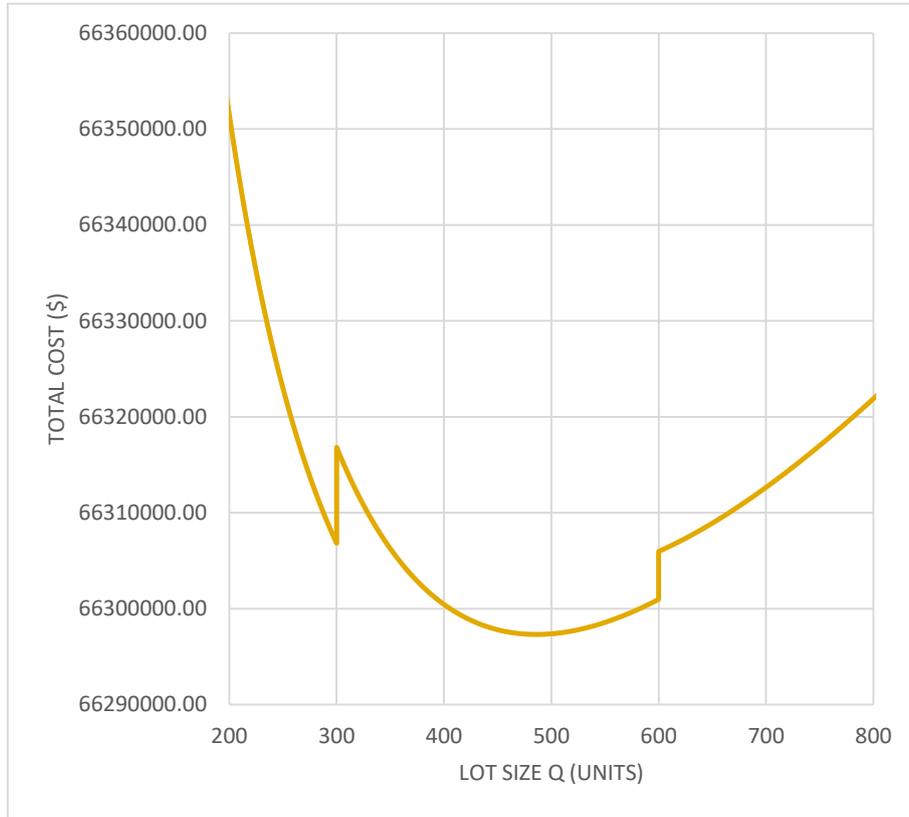

**Figure 2.** Total cost of the inventory system

To observe the variation in $Q^*_{opt}$ and $Q^*_{app}$, a sensibility analysis is carried out in the next section. Besides, it is also studied how the approximation is affected due to the change in other parameters and their effect in the total cost function.

### 4.3. Sensitivity analysis

This section presents a sensitivity analysis to observe the effect of the variation in the four parameters in two ways: the effect in the optimal lot size and in the accuracy of the approximation. As it is shown in Table 4, the percentage of variation in the value of each parameter is in the column $\Delta P(\%)$. Regarding the lot size, three calculations are made. First, the optimal lot size is given in $Q^*_{opt}$. Second, the variation in the optimal lot size is depicted in column $\Delta Q^*_{opt}(\%)$. Finally, the accuracy of the approximation in percentage is reported in the column Approx.(Q) which value is calculated as $\dfrac{Q^*_{opt} - Q^*_{app}}{Q^*_{opt}}$. The same calculations are made for the total cost function. The optimal total cost is presented in the column TOTAL COST (TC)*, its variation regarding the parameter's change in $\Delta TC(\%)$, and the accuracy of the corresponding approximation in column Approx.(C), where this latter value is calculated as $\dfrac{TC(Q^*_{opt}) - AC(Q^*_{app})}{TC(Q^*_{opt})}$

Table 4. Sensitivity analysis

| Parameter ($P$) | $\Delta P(\%)$ | $Q^*_{opt}$ | $\Delta Q^*_{opt}(\%)$ | Approx.(Q) (%) | TOTAL COST(TC)* | $\Delta TC(\%)$ | Approx.(C) (%) |
|---|---|---|---|---|---|---|---|
| $c$ | -20 | 486.0835 | 0.0000 | 0.0011 | 66272295.3469 | -0.0377 | 0.00000129 |
| | -10 | 486.0835 | 0.0000 | 0.0011 | 66284795.3469 | -0.0189 | 0.00000129 |
| | 0 | 486.0835 | 0.0000 | 0.0011 | 66297295.3469 | 0.0000 | 0.00000129 |
| | 10 | 486.0835 | 0.0000 | 0.0011 | 66309795.3469 | 0.0189 | 0.00000129 |
| | 20 | 486.0835 | 0.0000 | 0.0011 | 66322295.3469 | 0.0377 | 0.00000129 |
| $D$ | -20 | 434.7323 | -10.5643 | 0.0007 | 53053338.9666 | -19.9766 | 0.00000103 |
| | -10 | 461.1212 | -5.1354 | 0.0009 | 59675558.2774 | -9.9879 | 0.00000116 |
| | 0 | 486.0835 | 0.0000 | 0.0011 | 66297295.3469 | 0.0000 | 0.00000129 |
| | 10 | 509.8286 | 4.8850 | 0.0012 | 72918621.0249 | 9.9873 | 0.00000142 |
| | 20 | 532.5195 | 9.5531 | 0.0014 | 79539590.3523 | 19.9741 | 0.00000155 |
| $r$ | -20 | 486.0152 | -0.0141 | 0.000 | 66296672.7087 | -0.0009 | 0.00000066 |
| | -10 | 486.0474 | -0.0074 | 0.0008 | 66296983.3877 | -0.0005 | 0.00000094 |
| | 0 | 486.0835 | 0.0000 | 0.0011 | 66297295.3469 | 0.0000 | 0.00000129 |
| | 10 | 486.1235 | 0.0082 | 0.0014 | 66297608.5914 | 0.0005 | 0.00000172 |
| | 20 | 486.1675 | 0.0173 | 0.0018 | 66297923.1260 | 0.0009 | 0.00000223 |
| $l$ | -20 | 535.9360 | 10.2560 | 0.0008 | 66281389.6357 | -0.0240 | 0.00000085 |
| | -10 | 509.1865 | 4.7529 | 0.0009 | 66289528.7343 | -0.0117 | 0.00000106 |
| | 0 | 486.0835 | 0.0000 | 0.0011 | 66297295.3469 | 0.0000 | 0.00000129 |
| | 10 | 465.8674 | -4.1590 | 0.0012 | 66304738.1502 | 0.0112 | 0.00000155 |
| | 20 | 447.9831 | -7.8382 | 0.0014 | 66311896.0663 | 0.0220 | 0.00000182 |

Regarding the optimal lot size, it is hugely affected by the parameters $r$ and $D$. However, the approximation works in a very acceptable way. Regarding the optimal total cost, it is only vastly affected by the demand. However, the approximation works in a very satisfactory manner, again. Furthermore, three scenarios are analyzed in which the value of the parameter $r$ is changed to observe the behavior of the optimal lot size (now S-EOQ), the optimal lot size using equation (27) (now S-EOQ with approximation) and an optimal lot size if a company just takes into account environmental costs (now $Q_s$ from environmental costs).

**Scenario 1 ( $r = 0.004$ )**

In this scenario, a numerical experimentation is made with a low value for $r$. In Figure 3, the total cost curve and the environmental cost curve are plotted.

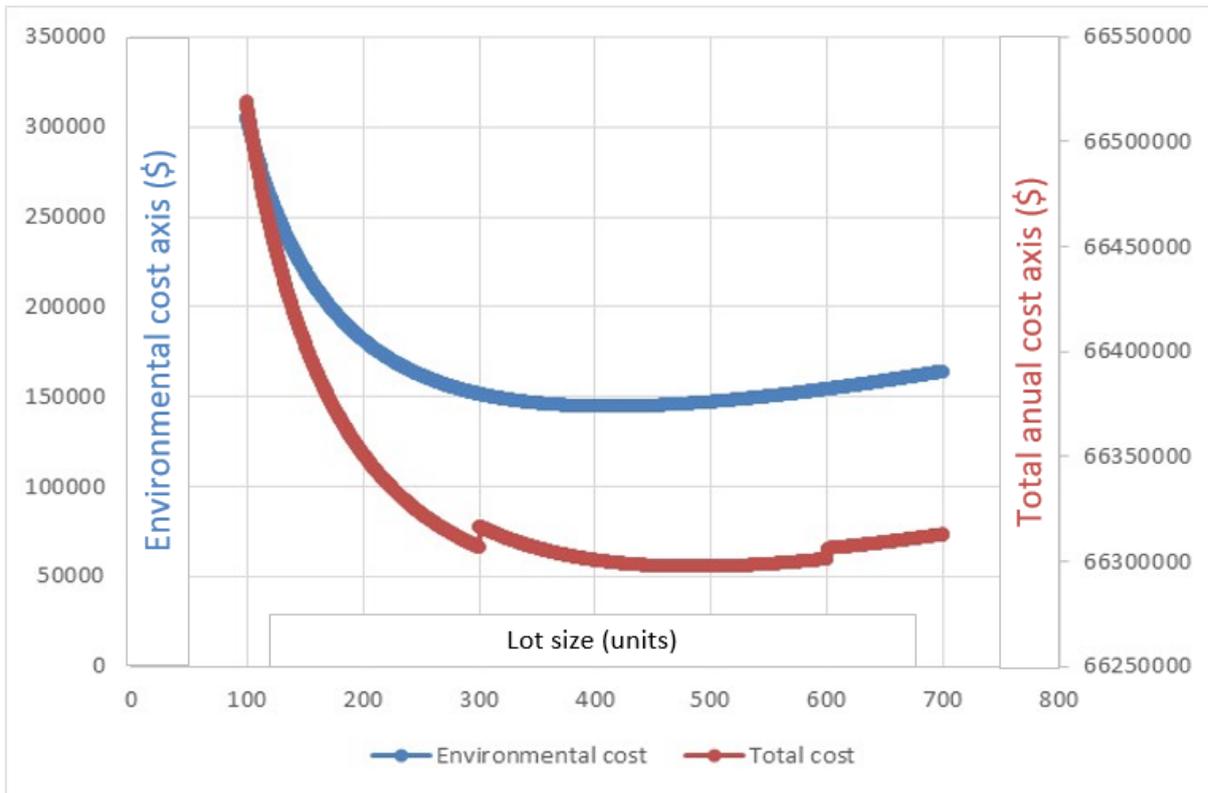

**Figure 3.** Total annual cost and environmental annual cost for Scenario 1

Table 5 shows the obtained results for this scenario. The difference between the sustainable lot size and the merely environmental lot size is 15.055%. It must be observed that the corresponding region in the total cost function in which the optimal lot size is found is slightly flat. Therefore, the variation in total cost is low.

**Table 5.** Results of scenario 1

| | |
|---|---:|
| S-EOQ | 486.084 |
| S-EOQ with approximation | 486.078 |
| $Q_s$ (from environmental costs) | 412.905 |

**Scenario 2 ( $r = 0.04$ )**

In this scenario, a numerical experimentation is made with a medium value for $r$. In Figure 4, the total and environmental cost curves are plotted.

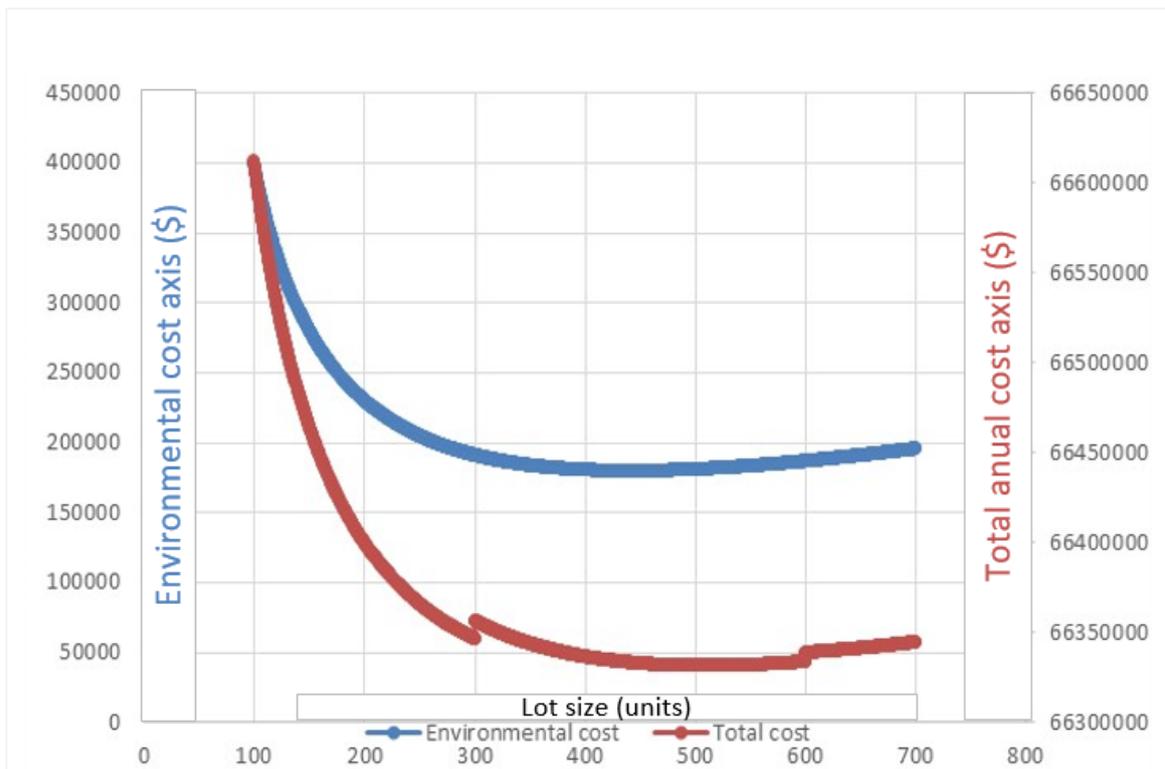

Figure 4. Total annual cost and environmental annual cost for Scenario 2

Table 6 shows the obtained results for this scenario. The difference between the sustainable lot size and the merely environmental lot size is 15.256%. It must be observed that the corresponding region in the total cost function is quite flat. Therefore, the variation in total costs is almost very low.

Table 6. Results of scenario 2

| | |
|---|---|
| S-EOQ | 509.002 |
| S-EOQ with approximation | 503.831 |
| $Q_s$ (from environmental costs) | 441.623 |

**Scenario 3 ($r = 0.2$)**

In this scenario, a numerical experimentation is made with a high value for $r$. In Figure 5, the total and the environmental cost functions are plotted.

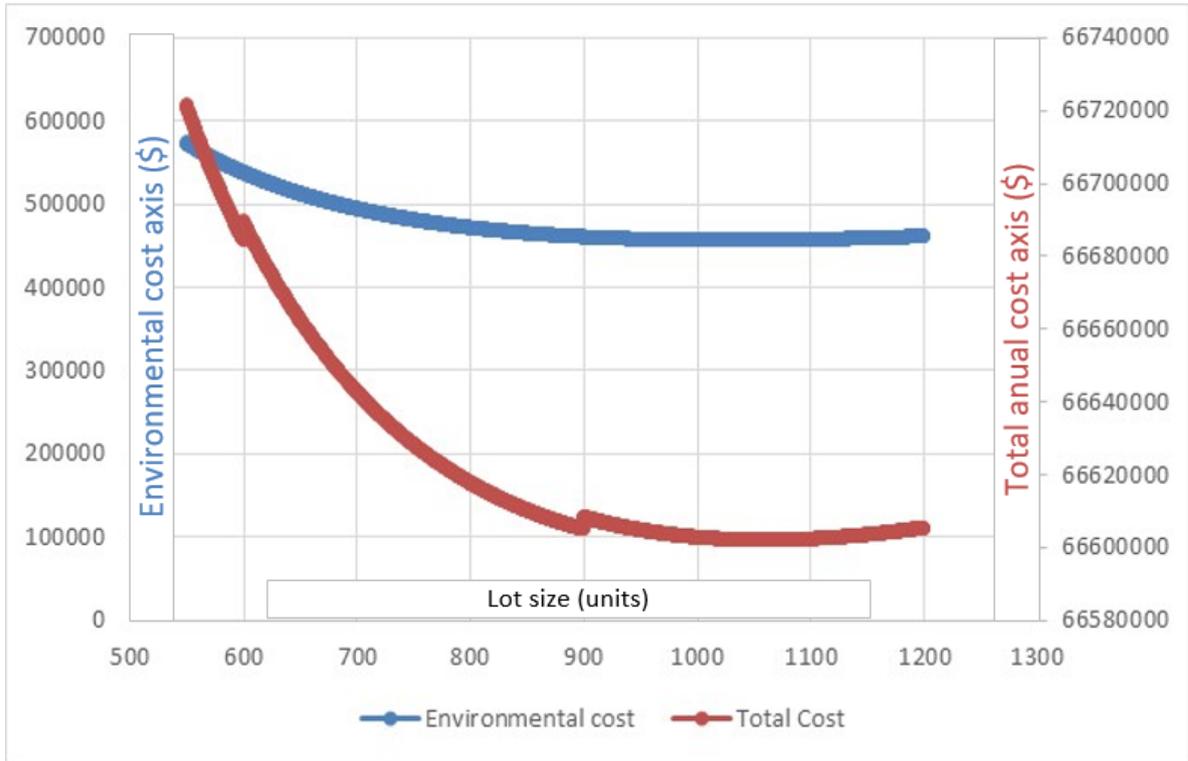

**Figure 5.** Total annual cost and environmental annual cost for Scenario 3

Table 7 summarizes the obtained results for this scenario. The difference between the sustainable lot size and the merely environmental lot size is 2.8301%. It is worth noting that the corresponding region in the total cost function is pretty flat. Therefore, the variation in total costs is negligible.

**Table 7**. Results of scenario 3

| | |
|---|---|
| S-EOQ | 1060 |
| S-EOQ with approximation | 1200 |
| $Q_s$ (from environmental costs) | 1030 |

It is observed that as the value of the parameter $r$ increases, the region in which the optimal lot size for the total cost has been found has become flatter. Then, the variation in the total cost function is negligible for high values of $r$, for any value of the lot size in the corresponding range. However, as it is concluded from the sensitivity analysis, as the parameter $r$ grows, the accuracy of the approximation diminishes. Then, a trade-off between the accuracy of the approximation and the possibility of obtaining a flat optimal-lot-size region needs to be considered. For that reason, in the following subsection, some observations about related to this topic are pointed.

### 4.4. Observations

According to numerical experiments, it is observed the following:

$rD$ must be less or equal than 78.8348 to be solved by Lingo. That means in our case ($D=5000$) that $r$ must be less or equal than 0.0157 to obtain a solution using the software.

$rD$ must be less or equal than 1000 to obtain a good enough approximation using a Taylor polynomial.

Due to these two empirical observations, the proposed methodology is a good approach to solve cases in which 78.8348< $rD$ <1000 for the selected parameters being intact.

As there is the necessity of finding integer lot sizes in order to be able of making an order in a real company, a formula to obtain an integer optimal solution is given in Section 5.

## 5. Case 2: when the lot size ($Q$) is an integer variable

In this section, the discretization algorithm to obtain the optimal solution for the problem stated in Section 3 is shown. Besides, a numerical experiment to compute the discrete optimal solution using the data in Table 1 is given.

### 5.1. Discretization algorithm

As it was proposed in García-Laguna et al. (2010), the discrete case in the EOQ inventory model is calculated by using marginal analysis in an easy-to-use approach. The use of that algorithm is preferred instead of the simple evaluation of the cost function at the two closest integer values to the continuous solution because it is more efficient, and it is used at each range $j = \{DP_j, DP_{j+1}\}$ in a similar way to the main procedure in Subsection 4.1. The solution procedure aims to solve the following optimization problem:

$$\text{Minimize } C(Q) = K\frac{D}{Q} + h\frac{Q}{2} + w$$

Subject to: $Q > 0$ and integer

In that sense, the first step to use the previous work is to convert the approximated total cost function, expressed in equation (25), in a similar way.

$$M(Q) = \left[(\varepsilon C_e + \frac{lr^2 D}{4}C_e + 2\beta\frac{d}{v} + y_0 + C_p DP_k + A + 2a)\frac{D}{Q}\right] + \left[(gC_e + lC_e + h)\frac{Q}{2}\right] + w \quad (31)$$

Where $w$ is given by:

$$w = \frac{l}{2}rDC_e + \gamma(\theta+\alpha)D + cD + bdD(1+\alpha)$$

Hence, the approximated cost function $M(Q)$ in equation (25) and (31) is expressed as:

$$M(Q) = K'\frac{D}{Q} + h'\frac{Q}{2} + w \tag{32}$$

Where $K'$ and $h'$ are given by:

$$K' = \varepsilon C_e + \frac{lr^2 D}{4}C_e + 2\beta\frac{d}{v} + y_0 + C_p DP_k + A + 2a$$
$$h' = gC_e + lC_e + h$$

Therefore, the expression is ready to apply in the following procedure:

**Step 1.** Set the $j$ th ranges and for each range $j = \{DP_j, DP_{j+1}\}$, calculate $Q'_j{}^*$ with following procedure:

    **Step 2.1** Compute $x = -0.5 + \sqrt{0.25 + 2K'D'/h'}$

    **Step 2.2** If $x$ is not an integer number, then $Q'_j{}^* = \lceil x \rceil$ is the unique optimal solution. Otherwise, both $Q'_j{}^* = \lceil x \rceil$ and $Q'_j{}^* = \lceil x \rceil + 1$ are optimal integer solutions.

    **Step 2.3** The minimum cost is obtained by calculating equation (32).

**Step 2.** For each range, the local optimal $Q_j{}^*$ is determined as follows:

$Q_j{}^* = Q'_j{}^*$ when $DP_j < Q'_j{}^* \leq DP_{j+1}$

If $Q'_j{}^*$ is not into the range $j = \{DP_j, DP_{j+1}\}$, then

$Q_j{}^* = DP_j$ when $C(DP_j) \leq C(DP_{j+1})$

$Q_j{}^* = DP_{j+1}$ when $C(DP_j) > C(DP_{j+1})$

**Step 3.** The total minimum cost is obtained as the minimum among the local total minimum cost values that satisfy $C(Q^*) = \min\{C(Q_j{}^*)\}\forall j$ and so $Q^*$ is the optimal value for the lot size.

In Subsection 5.2 a numerical example is carried out in order to clarify the use of the optimization algorithm.

### 5.2. Numerical experiment

This procedure is useful to obtain the results for the numerical example presented in this section. The only exception with the main procedure presented in Section 4 is that, since the lot size is now an integer variable, it is defined by using the lower case letter $q$ to emphasize the difference. The optimization procedure is reported in Table 8 and Table 9. The interpretation of these tables must be done according to the guidelines pointed for Table 2 and Table 3 previously discussed. Again, the iterative optimization procedure is carried out using the total and approximated cost function. The results for both options are reported in Table 8 and Table 9, respectively.

**Table 8**. Discretization algorithm -Total cost function

| $j+1$ | $n_1$ | $n_2$ | $DP_j$ | $DP_{j+1}$ | $TC(q'^*_j app)$ | $TC(DP_j)$ | $TC(DP_{j+1})$ | $q'^*_j app$ | Condition |
|---|---|---|---|---|---|---|---|---|---|
| 1 | 1 | 0 | 0 | 300 | | | 66306802.260 | 467 | F |
| 2 | 0 | 1 | 300 | 600 | **66297295.349** | | | **486** | T |
| 3 | 1 | 1 | 600 | 900 | | **66305950.560** | 66332800.248 | 504 | F |
| 4 | 2 | 1 | 900 | 1200 | | **66336133.582** | 66374075.139 | 521 | F |
| 5 | 1 | 2 | 1200 | 1500 | | **66376575.139** | 66419120.089 | 538 | F |
| 6 | 2 | 2 | 1500 | 1800 | | **66421120.089** | 66466050.062 | 554 | F |

**Table 9**. Discretization algorithm -approximated total cost function

| $j+1$ | $n_1$ | $n_2$ | $DP_j$ | $DP_{j+1}$ | $AC(q'^*_j app)$ | $AC(DP_j)$ | $AC(DP_{j+1})$ | $q'^*_j app$ | Condition |
|---|---|---|---|---|---|---|---|---|---|
| 1 | 1 | 0 | 0 | 300 | | | 66306800.000 | 467 | F |
| 2 | 0 | 1 | 300 | 600 | **66297294.494** | | | **486** | T |
| 3 | 1 | 1 | 600 | 900 | | **66305950.000** | 66332800.000 | 504 | F |
| 4 | 2 | 1 | 900 | 1200 | | **66336133.333** | 66374075.000 | 521 | F |
| 5 | 1 | 2 | 1200 | 1500 | | **66376575.000** | 66419120.000 | 538 | F |
| 6 | 2 | 2 | 1500 | 1800 | | **66421120.000** | 66466050.000 | 554 | F |

## 6. Conclusions

Global warming is a very important problem which is affecting the whole human race and its environment. For that reason, several global efforts aim at minimizing its impact by the use of agreements among different countries. One of the remarkable compromises stated in these documents is to develop strategies in different industrial activities to reduce the GHG emissions. Thus, sustainable inventory models represent remarkable academic contributions to reduce the environmental impact of the activities related to logistics management and to include the social and economic criteria into the decision making process.

The main contribution of this sustainable inventory model to the literature is the formulation of the emissions cost with an exponential relationship. However, the incorporated expression does not allow to obtain the optimal lot size in a close form formulae. For that reason, an approximation with a Taylor's polynomial is carried out to find an approximate total cost function and its corresponding optimal lot size. Then, the objective was to validate that the use of this approximation can be made without incurring in too big increment of costs.

The sustainable economic order quantity using the total cost function is 486.084 and by using the proposed approximation 486.078 (0.0011% of difference). Moreover, the difference between the total cost and the approximate cost is 0.00000129% which indicates that is a good enough approximation.

A sensitivity analysis was also carried out in which four parameters were studied:

a) $c$ does not affect neither the optimal quantity $Q^*$ nor any of the approximations. However, it is directly related with the total cost.

b) $D$ is directly related with the optimal quantity $Q^*$ and with the total cost. Both approximations, the sustainable economic order quantity and the total cost are better when $D$ is reduced.

c) $r$ is directly related with the optimal quantity $Q^*$ and with the total cost. Both approximations, the sustainable economic order quantity and the total cost are better when $r$ decreases.

d) $l$ is inversely related with the optimal quantity $Q^*$ and with the total cost. Both approximations, the sustainable economic order quantity and the total cost are better when $l$ is diminished.

During the numerical experimentation, three scenarios were studied, and it is possible to conclude that when $r$ increases the difference between $Q_s$ and the S-EOQ is reduced. On the one hand, it suggests that there are some cases in which it is possible to migrate to some completely sustainable police without incurring in high amounts of extra costs. On the other hand, it shows that there can be a trade-off between the sustainable benefit and the accuracy lost by increasing the parameter $r$.

Due to these two empirical observations, it was proved that the proposed methodology is a good approach to solve cases in which 78.8348< $rD$ <1000 considering the rest of the parameters intact. Another of the advantages of the proposed methodology is that Lingo is not able to find an optimal solution in this case. However, as it was depicted in Subsection 3.7, it is possible to find an approximation good enough.

For future work new carbon emissions functions can be developed and explored using Taylor series and some other approximation techniques. Moreover, some other characteristics such as imperfect quality and backorders can also be included in conjunction with the sustainability issues.

## Appendix A

**The carbon emissions function**

In this section, the parameters $r$ and $l$, which were introduced in Section 3, are explained and discussed. First, the new expression for the carbon emissions introduced in this work is shown and the benefit of using this new expression is commented. Then, mathematical relationships for $r$ and $l$ are developed. Finally, a practical methodology to understand and to obtain the parameter $r$ and $l$ is developed.

As it is mentioned in the mathematical formulation, there are some emissions $CO_2(Q)$ depending on the lot size that cannot be explained as the typical equation $CO_2(Q) = \varepsilon \frac{D}{Q} + g \frac{Q}{2}$ since the first term is decrescent and the last one is crescent regarding to the lot size and there are some activities or components of the inventory system that does not fit with this behavior. Then, a new expression is proposed as follows:

$$CO_2(Q) = \varepsilon \frac{D}{Q} + g \frac{Q}{2} + l \frac{Q}{2} e^{\frac{rD}{Q}}$$

The difference between the formulas is the addend $l\frac{Q}{2}e^{\frac{rD}{Q}}$. The amount of carbon emissions, according to this new expression, is supposed to grow while the lot size does tend to a straight line for lot sizes above a critical value that is introduced in this work as the critical lot size and represented as $Q^*$. It is also observed that the annual amount of emissions is increasing rapidly while the lot size decreases below $Q^*$.

The benefit of using this expression is the possibility of incorporating the emissions curve behavior of some machines that depends on the average inventory level when the lot size is above a certain value. However, when the number of orders during a period of time is above a critical number of orders $N^*$ and its corresponding lot size is below the critical lot size $Q^*$ the amount of emissions during the year increases rapidly in a very undesired way for the company purposes.

By modifying the parameters $r$ and $l$, the curve shape can fit into a variety of shapes for different company purposes. This new expression can model the behavior of some machines such as industrial air conditioners. If we assume that a plant needs to be turned on and down each cycle to make the replenishment of products, the number of orders affects considerably the environmental performance.

On the other hand, for large lot sizes, the amount of emissions during a period of time increases in a closed way to a straight line since the average inventory level is proportional to the amount of energy consumed by the air conditioner and its corresponding carbon footprint. However, when the number of cycles increases above its $N^*$, the environmental performance of the air conditioner is affected since it is turned on and down too many times. Then, the corresponding carbon emissions to that range of lot sizes is then increasing undesirably.

It is worth noting that fixing $r$ and $l$ the corresponding emissions curve could be modeled and the carbon emissions function proposed in Section 3 can be used with the algorithm introduced in this work (main procedure and discretization algorithm). Until now, the parameters $r$ and $l$ are introduced as shape parameters. However, it is useful to observe some mathematical relationships for both to improve the understanding of these parameters from a more practical point of view.

**Parameter $r$**

First, it is assumed that the emissions of some component is a function of the lot size and is $_fCO_2(Q)$. It must be observed that the corresponding carbon emissions curve has an absolute minimum in the critical lot size $Q^*$ as it was mentioned before. An expression for this value is obtained using the first derivative:

$$_fCO_2(Q) = l\frac{Q}{2}e^{\frac{rD}{Q}}$$

$$_fCO'_2(Q) = \frac{l}{2}e^{\frac{rD}{Q}} - \frac{l}{2}Qe^{\frac{rD}{Q}}\left(\frac{rD}{Q^2}\right)$$

$$_fCO'_2(Q) = \frac{l}{2}e^{\frac{rD}{Q}}\left(1 - \frac{rD}{Q}\right)$$

Then, making $_fCO'_2(Q^*) = 0$:

$$1 - \frac{rD}{Q^*} = 0$$

$$rD = Q^*$$

Finally, $r$ is calculated as $r = \frac{Q^*}{D}$ which is the corresponding length cycle to the critical lot size.

**Parameter $l$**

First, the behavior of the graph when $Q > Q^*$ is analyzed. It is observed that the curve tends to be a straight line with slope $m$

$$m = \lim_{Q \to \infty}\left(\frac{_fCO_2(Q)}{Q}\right)$$

$$m = \lim_{Q \to \infty}\left(\frac{l\frac{Q}{2}e^{\frac{rD}{Q}}}{Q}\right)$$

$$m = \left(\frac{l}{2}e^0\right) = \frac{l}{2}$$

Then, the asymptote of the curve has a slope $m = \frac{l}{2}$. Thus, it is stated that when $Q \to \infty$, the rate of change of the annual emissions is with respect to the lot size is:

$$\frac{\Delta_f CO_2(Q)}{\Delta Q} = \frac{l}{2}$$

And it can be rewritten as:

$$\frac{\Delta_f CO_2(Q)}{\Delta \frac{Q}{2}} = l$$

Since $\frac{Q}{2}$ is the average inventory level, $l$ represents the rate of change of the annual emissions of the component with respect to the average inventory level.

**Practical methodology**

A practical way of determining the parameters $r$ and $l$ for the emissions of some equipment (as a machine), assuming that it can be modeled as $_fCO_2(Q)$ is as follows.

1. Determine the critical number of orders $N^*$ above which the environmental performance starts growing rapidly. Then, calculate the critical lot size $Q^*$ and the parameter $r$ using $Q^* = \frac{D}{N^*}$ and $r = \frac{Q^*}{D}$.
2. Determine the parameter $l$ as the rate of change of the carbon emissions with respect to the average inventory level for the range in which the lot sizes are very much higher than the critical lot size $Q^*$.
3. With the parameters $r$ and $l$ $_fCO_2(Q)$ is defined. The next step is to verify that this emissions curve is adequate to model the real performance of the equipment.

**Example**

Let suppose that an engineer observes that an air-conditioner carbon footprint is proportional to the lot size for big lot sizes. However, when the number of orders during the year increases above 100 orders, the annual emissions starts increasing and becoming undesirable so rapidly. The demand is 1000 units/year and it is known that above a lot size of 40 units, the annual emissions increases in 2 kg of $CO_2$ for each additional unit in the average inventory level.

**Solution**

The parameter $r$ and the critical lot size $Q^*$ are calculated:

$Q^* = $ 1000(units/year)/100(orders/year) = 10 (units/order)

$r = $ 10(units/order)/1000(units/year) = 0.01 (year/order)

The parameter $l$ is the rate of change of the emissions with respect to the average inventory level. For this example, it is assumed that 40 is very much higher than 10. Then:

$l = $ 2 kg of $CO_2$/year/unit

Then, the function $_fCO_2(Q) = l\frac{Q}{2}e^{\frac{rD}{Q}}$ is $_fCO_2(Q) = 2\frac{Q}{2}e^{\frac{0.01(1000)}{Q}}$ and finally:

$$_fCO_2(Q) = Qe^{\frac{10}{Q}} \text{ (kg of } CO_2/\text{year)}$$

The last step is to verify that $_fCO_2(Q)$ models the emissions of the air conditioner. The last example includes the assumption of that the emissions curve tends to the asymptote close enough by above 40 units. However, it is possible to determine how far from the critical lot size the analysis need to be made. Defining that the curve emissions is close enough to the asymptote when $\frac{_fCO'_2(Q)}{m_a} \geq 0.9$, where $m_a$ is the asymptote slope, the following relationship is given:

$$0.9\frac{l}{2}Q \leq \frac{l}{2}e^{\frac{rD}{Q}}(Q - rD)$$

Since $rD = Q^*$, a replacement is made:

$$0.9Q \leq e^{\frac{Q^*}{Q}}(Q - Q^*)$$

It is observed that the last condition is satisfied when $Q \geq 3Q^*$.

**Acknowledgements**

To people around the world who are fighting the coronavirus using their best.